\documentclass[letterpaper, 10 pt, conference]{ieeeconf}  

\IEEEoverridecommandlockouts                              

\newtheorem{innercustomthm}{Theorem}
\newenvironment{customthm}[1]
  {\renewcommand\theinnercustomthm{#1}\innercustomthm}
  {\endinnercustomthm}

\newtheorem{innercustomass}{Assumption}
\newenvironment{customass}[1]
  {\renewcommand\theinnercustomass{#1}\innercustomass}
  {\endinnercustomass}  
  
\newtheorem{innercustomcor}{Corollary}
\newenvironment{customcor}[1]
  {\renewcommand\theinnercustomcor{#1}\innercustomcor}
  {\endinnercustomcor}

\usepackage{graphics} 
\usepackage{epsfig} 
\usepackage{cite}
\usepackage{multirow}
\usepackage{comment}

\usepackage{amsmath} 
\usepackage{amssymb}  

\usepackage{mathrsfs}
\usepackage{todonotes}

\title{\LARGE \bf
Conditions for State and Control Constraint Activation in Coordination of Connected and Automated Vehicles}

\author{A M Ishtiaque Mahbub, {\itshape{Student Member, IEEE}}, Andreas A. Malikopoulos, {\itshape{Senior Member, IEEE}}
\thanks{This research was supported by ARPAE's NEXTCAR program under the award number DE-AR0000796.}
\thanks{The authors are with the Department of Mechanical Engineering, University of Delaware, Newark, DE 19716 USA (emails: \tt\small{mahbub@udel.edu};\tt\small{andreas@udel.edu}.)}}

\begin{document}

\maketitle
\thispagestyle{empty}
\pagenumbering{arabic}

\begin{abstract}
Connected and automated vehicles (CAVs) provide the most intriguing opportunity to reduce pollution, energy consumption, and travel delays. In earlier work, we addressed the optimal coordination of CAVs using Hamiltonian analysis. In this paper, we investigate the nature of the unconstrained problem and provide conditions under which the state and control constraints become active. We derive a closed-form analytical solution of the constrained optimization problem and evaluate the solution using numerical simulation.

\end{abstract}

\indent


\section{Introduction} \label{sec:1}
The implementation of an emerging transportation system with connected automated vehicles (CAVs) enables a novel computational framework to provide real-time control actions that optimizes energy consumption and associated benefits. CAVs can alleviate congestion at major transportation segments \cite{Margiotta2011, mahbub2020sae-1}, reduce emission, improve fuel efficiency \cite{Lee2013, Malikopoulos2017}, and increase passenger safety. 
Several research efforts have been reported in the literature proposing either centralized or decentralized approaches on coordinating CAVs at different traffic scenarios, e.g., intersections, merging at roadways and roundabouts, speed reduction zones, to improve traffic flow and reduce stop-and-go driving. Some papers have used a reservation scheme \cite{Dresner2008, DeLaFortelle2010, Huang2012} to control a signal-free intersection of two roads. Other approaches have focused on coordinating vehicles to improve the traffic flow \cite{Zohdy2012, Yan2009, kim2014}. 

More recently, an optimal control framework was established for coordinating online CAVs. A closed-form, analytical solution without considering state and control constraints was presented in \cite{Ntousakis:2016aa} for coordinating online CAVs at highway on-ramps, in \cite{Mahbub2019ACC} at intersections, and in \cite{Zhao2019CCTA-1} for traffic corridors. The solution to the state and control unconstrained optimization problem discussed above shows possible state and control constraint activation within the optimization horizon and at the boundaries. To mitigate terminal jerk, Ntousakis et al. \cite{Ntousakis:2016aa} reformulated the optimal control problem with terminal constraint on state and control, which does not guarantee constrained state/control values within the optimization horizon. The optimal control problem considering state and control constraints was addressed in \cite{Malikopoulos2017} at an urban intersection, where the resulting formulation relies on piecing the unconstrained and constrained arcs together, and leads to a recursive computation of the solution until all constraint activation cases are resolved. This procedure may end up in a complicated recursive structure, and can become computationally exhaustive due to the presence of implicit system of equations to be solved numerically, preventing possible real-time implementation. The process is similar if the rear-end safety constraint is included \cite{malikopoulos2019ACC}.

In this paper, we investigate the formulation presented in \cite{Malikopoulos2017} to provide useful insights about the state and control constraint activation to increase computational efficiency, and derive closed form analytical solution using the Hamiltonian analysis.
In particular, we provide (1) conditions under which the constraint activation space can be reduced, and (2) conditions under which specific combination of state and control constraint activation can be identified to increase computational efficiency.

The remainder of the paper is organized as follows. In Section II, we introduce the problem formulation and present the unconstrained case. In Section III, we discuss different aspects of the state and control constrained formulation in detail. In Section IV, we provide the analytical solution of the constrained optimization. In Section V, we evaluate the effectiveness of the proposed approach in a simulation environment. We draw conclusions and discuss next steps in Section VI.

\section{Problem Formulation} \label{sec:2}

\begin{figure}[ht]
\centering
\includegraphics[scale=0.35]{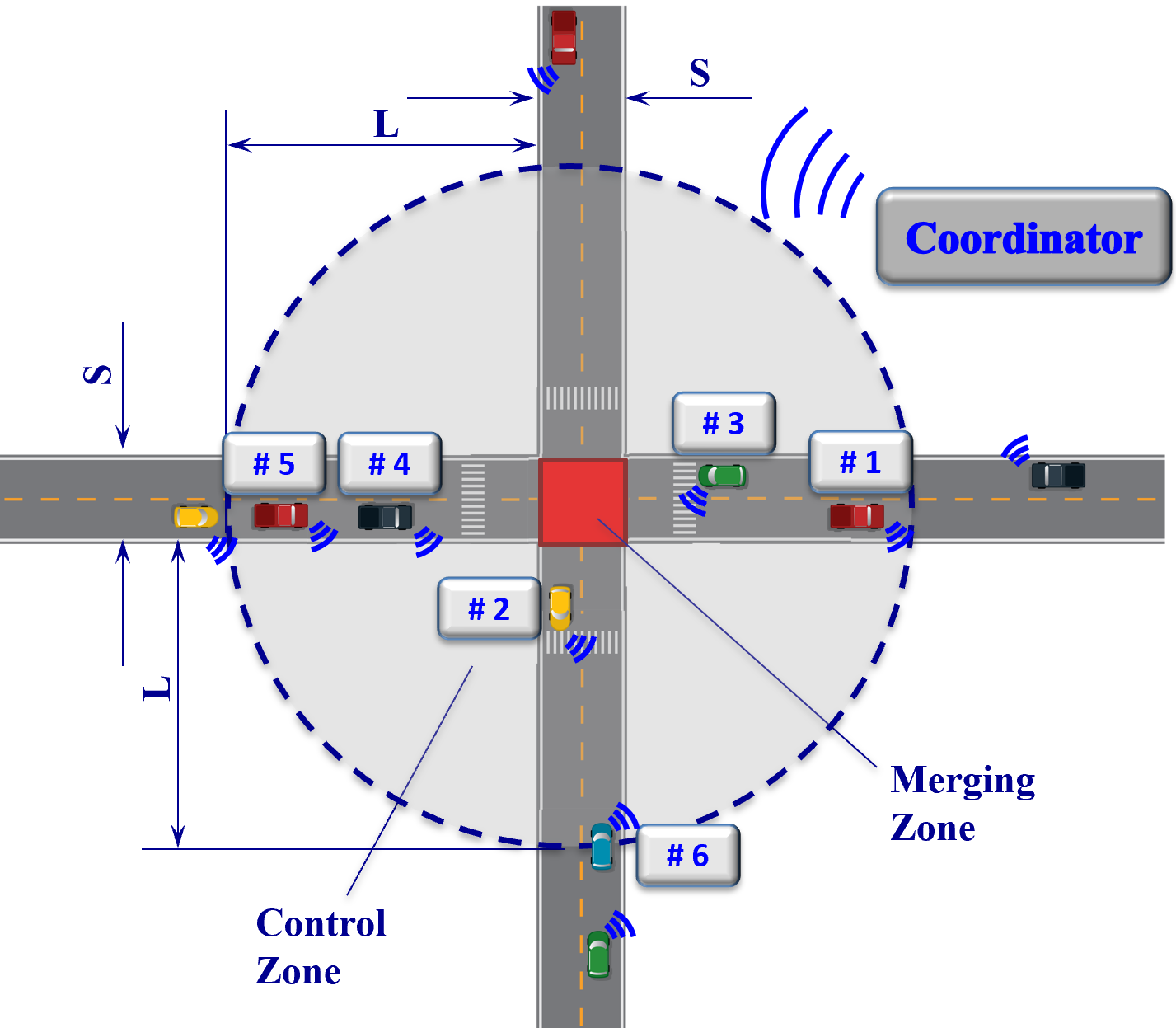} \caption{Connected automated vehicles at a four-way intersection.}%
\label{fig:4}%
\end{figure}

Although our theoretical framework can be applied to any
traffic scenario, e.g., merging at roadways \cite{Ntousakis:2016aa} and roundabouts, passing through speed reduction
zones \cite{Malikopoulos2018c}, we use an
intersection (Fig. \ref{fig:4}) as a reference to present the fundamental ideas
and results of this paper, since an intersection provides unique features
making it technically more challenging compared to
other traffic scenarios.

We define the area illustrated by the red square of dimension $S$ in Fig. \ref{fig:4} as the \textit{merging zone} where possible lateral collision may occur. Upstream of the merging zone, we define a \textit{control zone} of length $L$, inside of which CAVs can coordinate with each other before they pass the merging zone. The intersection has a coordinator that only maintains the queue of the CAVs inside the control zone, and is not involved in any decision-making process. When a CAV enters the control zone, the coordinator receives its information and assigns a unique identity $i\in\mathbb{N}$ to it. 
The objective of each CAV is to derive its optimal control input to cross the intersection without activating any of the safety, state and control constraints.

\subsection{Vehicle Dynamics and Constraints}
Let $\mathcal{N}(t)=\{1,..., N(t)\}$, where $N(t)\in\mathbb{N}$ is the number of CAVs inside the control zone at time $t\in\mathbb{R}^{+}$, be the queue of CAVs inside the control zone.
We model each CAV $i\in\mathcal{N}(t)$ as a double integrator
\begin{equation}%
\begin{split}
\dot{p}_{i}  =v_{i}(t), ~ \dot{v}_{i}   =u_{i}(t),~ \forall t\in [t_i^0, t_i^f],\label{eq:model2}
\end{split}
\end{equation}
where $p_{i}(t)\in\mathcal{P}_{i}$, $v_{i}(t)\in\mathcal{V}_{i}$, and $u_{i}(t)\in\mathcal{U}_{i}$ denote the position, speed and
acceleration (control input) of each CAV $i$ in the corridor. 
The sets $\mathcal{P}_{i}$, $\mathcal{V}_{i}$, and $\mathcal{U}_{i}$, $i\in\mathcal{N}(t),$ are complete and totally bounded subsets of $\mathbb{R}$.
Let $\textbf{x}_{i}(t)=\left[p_{i}(t) ~ v_{i}(t)\right] ^{T}$ denote the state vector of each CAV $i\in \mathcal{N}(t)$, with initial value $\textbf{x}_{i}^{0}=\left[p_{i}^{0} ~ v_{i}^{0}\right]  ^{T}$ taking values in $\mathcal{X}_{i}%
=\mathcal{P}_{i}\times\mathcal{V}_{i}$. The state space $\mathcal{X}_{i}$ for each CAV $i\in \mathcal{N}(t)$ is closed with respect to the induced topology on $\mathcal{P}_{i}\times \mathcal{V}_{i}$ and thus, it is compact. Each CAV $i$ enters the control zone at $t_i^0$, enters the merging zone at $t_i^{m}$, and exits the merging zone at $t_i^{f}$.

To ensure that the control input and vehicle speed are within a
given admissible range, the following constraints are imposed
\begin{gather}%
u_{min} \leq u_{i}(t)\leq u_{max},\quad\text{and} \nonumber \\ 
0  \leq v_{min}\leq v_{i}(t)\leq v_{max},\quad\forall t\in\lbrack t_{i}%
^{0},t_{i}^{f}],
\label{eq:state_control_constraint}%
\end{gather}
where $u_{min}$, $u_{max}$ are the minimum and maximum
acceleration for each CAV $i\in\mathcal{N}(t)$, and $v_{min}$, $v_{max}$ are the minimum and maximum speed limits respectively.
 
We impose the following rear-end safety constraint 
\begin{equation}
s_{i}(t)=p_{k}(t)-p_{i}(t) \geq \delta_i(t),~ \forall t\in [t_i^0, t_i^f],
\label{eq:rearend_constraint}
\end{equation}
where $s_{i}(t)$ denotes the distance of CAV $i$ from CAV $k$ which is physically immediately ahead of $i$, and  $\delta_i(t)$ denotes minimum safe distance which is a function of speed $v_i(t)$.

In the modeling framework described above, we impose the following assumptions:\newline
\textbf{Assumption 1.}
Each CAV is equipped with sensors to measure and share their local information.\newline
\textbf{Assumption 2.}
Communication among CAVs occurs without any delays or errors. 

Assumption 1 restricts the problem to the case of 100\% CAV penetration rate.
Assumption 2 may be strong, but it is relatively straightforward to relax as long as the noise in the measurements and/or delays is bounded. For example, we can determine upper bounds on the state uncertainties as a result of sensing or communication errors and delays, and incorporate these into more conservative safety constraints.

\subsection{Upper-Level Vehicle Coordination}
The queue of CAVs inside the control zone $\mathcal{N}(t)$ can be determined as an outcome of an upper-level vehicle coordination problem as described in \cite{Malikopoulos2019CDC,Mahbub2019ACC}. 
In what follows, we assume that the time $t_i^m$ that each CAV enters the merging zone is determined as the solution of some upper-level coordination problem, and will focus only on a low-level energy optimization problem that will yield the optimal control input $u_i(t)$ for each CAV $i\in\mathcal{N}(t)$ to achieve the assigned merging time $t_i^m$.

\subsection{Low-Level Energy Optimization Problem}
For each CAV $i\in\mathcal{N}(t)$ traveling inside the control zone, we formulate an optimal control problem to minimize the $L^2$-norm of the control input, i.e., $\frac{1}{2}u_i^2(t)$, so that the transient engine operation can be minimized, and thus, fuel consumption  \cite{Malikopoulos2017},

\begin{gather}\label{eq:decentral_problem}
\min_{u_i\in U_i}J_{i}(u(t))=  \int_{t_i^{0}}^{t_i^{m}}\frac{1}{2}u_i^2(t)~dt,\\ 
\text{subject to}%
:\eqref{eq:model2},\eqref{eq:state_control_constraint},\nonumber\\
\text{and given }t_i^{0}\text{, }\text{
}p_{i}(t_i^{0})=0, v_{i}(t_i^{0})\text{, }t_i^{m},~ p_{i}(t_i^{m}).\nonumber
\end{gather}
Note that, we do not provide a desired speed $v_{i}(t_i^{m})$ at $t_i^m$. Additionally, we do not explicitly include the lateral and rear-end \eqref{eq:rearend_constraint} safety constraints. We enforce the lateral collision constraint by selecting the appropriate merging time $t_i^m$ for each CAV $i$ in the upper-level vehicle coordination problem. The activation of rear-end safety constraint can be avoided under proper initial conditions $[t_i^0,v_i(t_i^0)]$ as it was shown in \cite{Malikopoulos2018c}. 
From \eqref{eq:decentral_problem} and the state equations \eqref{eq:model2}, we formulate the state and control constraint \eqref{eq:state_control_constraint} adjoined Hamiltonian function $H_{i}\big(t, \boldsymbol{\lambda}(t),\boldsymbol{\mu}(t), \textbf{x}(t),u(t)\big)$ for each CAV $i\in\mathcal{N}(t)$,
\begin{gather}
H_{i}\big(t, \boldsymbol{\lambda}(t),\boldsymbol{\mu}(t), \textbf{x}(t),u(t)\big) 
=\frac{1}{2} u^{2}_{i} \nonumber
+ \lambda^{p}_{i} \cdot
v_{i} + \lambda^{v}_{i} \cdot u_{i} \\ \nonumber
+\mu^{a}_{i} \cdot(u_{i} - u_{max})+ \mu^{b}_{i} \cdot(u_{min} - u_{i})\\
+ \mu^{c}_{i} \cdot(v_{i} - v_{max}) 
+ \mu^{d}_{i} \cdot(v_{min} - v_{i}),
\quad \forall t \in[t_i^0, t_i^m], \label{eq:lagrangian}%
\end{gather}
where $\lambda^{p}_{i},~  \lambda^{v}_{i}$ are the co-state components corresponding to the state vector $\textbf{x}_i(t)$,
and $\mu^{a}_{i}$, $\mu^{b}_{i}$, $\mu^{c}_{i}$, $\mu^{d}_{i}$ are the Lagrange multipliers with the following KKT conditions, $\mu^{a}_{i} \cdot (u_{i}(t) - u_{max}) = 0$, $\mu^{b}_{i} \cdot (u_{min} - u_i(t)) = 0$, $\mu^{c}_{i} \cdot (v_{i}(t) - v_{max}) = 0$, $\mu^{d}_{i} \cdot (v_{min} - v_i(t)) = 0$, and $\mu^{a}_{i}$, $\mu^{b}_{i}$, $\mu^{c}_{i}$, $\mu^{d}_{i}>0$ on the constrained arc. 
The Euler-Lagrange equations are
\begin{gather}\label{eq:EL1}
\dot\lambda^{p}_{i}(t) = - \frac{\partial H_i}{\partial p_{i}} = 0, ~
\dot\lambda^{v}_{i} = - \frac{\partial H_i}{\partial v_{i}} =-\lambda^{p}_{i}-\mu^{c}_{i}+\mu^{d}_{i},\\
\frac{\partial H_i}{\partial u_{i}} = u_{i} + \lambda
^{v}_{i} + \mu^{a}_{i} - \mu^{b}_{i} = 0.\label{eq:nece_cond}
\end{gather}

\subsection{Unconstrained Optimization}
If the inequality state and control constraints \eqref{eq:state_control_constraint} are not active, we have $\mu^{a}_{i} = \mu^{b}_{i}= \mu^{c}_{i}=\mu^{d}_{i}=0$.
From (\ref{eq:EL1}) and \eqref{eq:nece_cond}, for each CAV $i\in \mathcal{N}(t)$ we derive the optimal control input $u_i^*(t)$, and the optimal states $p_i^*(t),~v_i^*(t)$ as
\begin{gather}
u^{*}_{i}(t) = a_{i} \cdot t + b_{i} ,~ \forall t\in [t_i^0, t_i^m], \label{eq:20}\\
v^{*}_{i}(t) = \frac{1}{2} a_{i} \cdot t^2 + b_{i} \cdot t + c_{i},~ \forall t\in [t_i^0, t_i^m],\label{eq:21}\\
p^{*}_{i}(t) = \frac{1}{6}  a_{i} \cdot t^3 +\frac{1}{2} b_{i} \cdot t^2 + c_{i}\cdot t + d_{i}, ~ \forall t\in [t_i^0, t_i^m], \label{eq:22}%
\end{gather}
where $a_i$, $b_{i}$, $c_{i}$, and $d_{i}$ are constants of integration which can be computed by using the boundary conditions in \eqref{eq:decentral_problem}.

\section{Constrained Optimization}
To derive a closed-form analytical solution for \eqref{eq:lagrangian}, we (1) identify the conditions for constraint(s) exclusion, (2) define the condition under which they become active, and (3) derive the final constrained optimal solution without recursion. 

\subsection{Condition of Constraint Exclusion}\label{sec:step1}
Although we have two state and two control constraints from \eqref{eq:state_control_constraint}, there are 15 constraint combinations in total that can be activated. In this section, we show that it is only possible for a subset of the constraints to be active within an unconstrained solution. Therefore, it is not necessary to consider all the cases in the constrained problem. In what follows, we delve deeper into the nature of the uncontrolled optimal solution to derive useful information about the possible existence of constraint activation within the control zone.

\textbf{Lemma 1.}
If $v_i(t_i^m)$ is not prescribed, and the terminal time $t_i^m$ is fixed, then
\begin{equation}
    t_i^m=-\frac{b_i}{a_i}. \label{eq:b_i}
\end{equation}

\begin{proof}
For fixed terminal time $t_i^m$ and unspecified state $v_i(t_i^m)$, $\lambda_i^v(t_i^m) = 0$.
Using $\lambda_i^v(t_i^m)$ in \eqref{eq:20}, we have $u_i(t_i^m)=a_it_i^m+b_i=0$,
and the result follows.
\end{proof} 

\textbf{Corollary 1.}
The constants $a_i$ and $b_i$ always have opposite signs.

\textbf{Lemma 2.}
The unconstrained optimal control input $u_i(t)$ is linearly decreasing if $v_i(t_i^0)<\frac{(p_i(t_i^m)-p_i(t_i^0))}{t_i^m}$.
The unconstrained optimal control input $u_i(t)$ is linearly increasing if $v_i(t_i^0)>\frac{(p_i(t_i^m)-p_i(t_i^0))}{t_i^m}$.

\begin{proof}
From \eqref{eq:21} and \eqref{eq:22}, we can write $v_i(t_i^0)=\frac{1}{2}a_i\cdot(t_i^0)^{2}+b_i\cdot t_i^0+c_i$ and $p_{i}(t_i^0) = \frac{1}{6}  a_{i} \cdot (t_i^0)^3 +\frac{1}{2} b_{i} \cdot (t_i^0)^2 + c_{i}\cdot t_i^0 + d_{i}$.
With $t_i^0=0$, we have
\begin{gather}
    c_i=v_i(t_i^0),~ d_i = p_i(t_i^0). \label{eq:c_i}
\end{gather}
From \eqref{eq:22}, we have $p_i(t_i^m)=\frac{1}{6}a_i\cdot(t_i^m)^3+\frac{1}{2}b_i\cdot(t_i^m)^2+c_i\cdot t_i^m+d_i.$
Using the results from \eqref{eq:c_i}, and \eqref{eq:b_i} in the above equation and solving for $a_i$, we have
\begin{equation}
  a_i=\frac{3(v_i(t_i^0)t_i^m-(p_i(t_i^m)-p_i(t_i^0)))}{(t_i^m)^3}. \label{eq:a_i}
\end{equation}
Since $t_i^m>0$, $a_i$ is non positive if $v_i(t_i^0)t_i^m-(p_i(t_i^m)-p_i(t_i^0))<0$,
which in turn implies negative slope of the linear control input $u_i^*(t)$, i.e., linearly decreasing acceleration profile. The second part of Lemma 2 can be proved following similar steps.
\end{proof} 

\textbf{Lemma 3.}
For a given unconstrained optimal profile, only a subset of state-control constraints can remain active. Furthermore, if either $v_i(t)-v_{max}\le0$ or $u_i(t)-u_{max}\le0$ is activated, none of the constraint pairs $v_{min}-v_{i}(t)\le0$ and $u_{min}-u_{i}(t)\le0$ can be activated. The reverse also holds.

\begin{proof} 
The optimal solution derived from the unconstrained case yields a linear control profile increasing or decreasing to the terminal value zero (Lemma 2). If the acceleration profile is linearly decreasing, we have $u_i(t)=a_it+b_i\ge 0>u_{min}$, $\forall t \in [t_i^0,t_i^m]$, implying that the constraint $u_{min}-u_i(t)\le0$ will never be activated. Now, applying the necessary condition of optimality in \eqref{eq:21}, we have $\frac{\partial v_i(t)}{\partial t} = a_i t + b_i=0$, which yields the inflection point at time $t=-\frac{b_i}{a_i}$ corresponding to the vertex of the concave parabola of \eqref{eq:21}, and equal to $t_i^m$ according to Lemma 1. As the inflection point of the concave parabola is located at $t_i^m$ and $v_{min}<v_i^0<v_{max}$, we have $v_i(t)>v_{min}$, $\forall t \in [t_i^0,t_i^m]$. This concludes the proof of the first part of Lemma 3. The second part of Lemma 3 can be proved following similar steps.
\end{proof}

\textbf{Remark 1.}
The sign of the constant $a_i$ provides insight for the possible state and control constraint activation.

The following result provides the condition under which certain state and control constraints never become active.

\textbf{Theorem 1.}
(i) The state constraint $v_{min}-v_{i}(t)\le 0$ and the control constraint $u_{min}-u_{i}(t)\le 0$ are not activated if $v_i(t_i^0)<\frac{(p_i(t_i^m)-p_i(t_i^0))}{t_i^m}$. (ii) The state constraint $v_i(t)-v_{max}\le 0$ and the control constraint $u_i(t)-u_{max}\le 0$ are not activated, if $v_i(t_i^0)>\frac{(p_i(t_i^m)-p_i(t_i^0))}{t_i^m}$.

\begin{proof}
If $v_i(t_i^0)<\frac{(p_i(t_i^m)-p_i(t_i^0))}{t_i^m}$ holds, then from Lemma 2, $a_i$ is negative and the optimal control input $u_i(t)$ is linearly decreasing. From Lemma 3, a decreasing optimal control input $u_i(t)$ indicates that the state constraint $v_{min}-v_{i}(t)\le 0$ and the control constraint $u_{min}-u_{i}(t)\le 0$ will never be activated, which concludes the proof of the first part. The second part of Theorem 1 can be proved following similar procedure.
\end{proof}

\subsection{Conditions of Constraint Activation}\label{sec:step2}
The following results provide the condition for which specific state and control constraints are activated.

\textbf{Lemma 4.}
Activation of the control constraint $u_i(t)-u_{max}\le 0$ or $u_{min}-u_{i}(t)\le 0$ occurs at time $t=t_i^0$.

\begin{proof}
For $a_i<0$, there is a possibility that either the state constraint $v_i(t)-v_{max}\le 0$, or the control constraint $u_i(t)-u_{max}\le 0$, or both to be activated (Lemma 3). In this case, the control input linearly decreases to zero at time $t=t_i^m$. Therefore, any activation of the control constraint $u_i(t)-u_{max}\le 0$ occurs at time $t=t_i^0$. Following similar reasoning, it can be proved that the activation of $u_{min}-u_{i}(t)\le 0$ occurs at time $t=t_i^0$
\end{proof}

\textbf{Theorem 2.}
(i) If $a_i<0$, the state constraint  $v_i(t)-v_{max}\le 0$ is activated if $t_i^m<\frac{3(p_i(t_i^m)-p_i(t_i^0))}{v_i^0+2v_{max}}$.
(ii) If $a_i>0$, the state constraint $v_{min}-v_{i}(t)\le 0$ is activated if $t_i^m>\frac{3(p_i(t_i^m)-p_i(t_i^0))}{v_i^0+2v_{min}}$.

\begin{proof}
Suppose that for $a_i<0$, there exists a time $t_s \in (t_i^0,t_i^m]$ at which the state constraint $v_i(t)-v_{max}\le 0$ is activated. Using \eqref{eq:c_i}, from \eqref{eq:21} we have $\frac{1}{2}a_i t_s^2+b_i t_s+v_i^0=v_{max}$.
Solving for $t_s$, we have $t_s = \frac{-2b_i\pm \sqrt{4b_i^2-8a_i\cdot(v_i^0-v_{max})}}{2a_i}$, yielding two possible solutions, either $t_{s,1}=t_i^m+ \sqrt{\frac{4b_i^2-8a_i\cdot(v_i^0-v_{max})}{4a_i^2}}$ or $t_{s,2}=t_i^m-\sqrt{\frac{4b_i^2-8a_i\cdot(v_i^0-v_{max})}{4a_i^2}}$.
Since $t_{s,1}$ is not feasible as $\sqrt{4b_i^2-8a_i\cdot(v_i^0-v_{max})}$ cannot be negative or equal to zero to satisfy $t_s<t_i^m$, $t_{s,2}$ is the only acceptable solution. To satisfy $t_s<t_i^m$, we have $\sqrt{4b_i^2-8a_i\cdot(v_i^0-v_{max})}>0$ resulting in $a_i<\frac{2(v_i^0-v_{max})}{(t_i^m)^2}$.
Using \eqref{eq:a_i}, the first part of Theorem 2 follows. The second part of Theorem 2 can be proved for $a_i>0$ following similar arguments.
\end{proof} 

\textbf{Theorem 3.}
(i) For $a_i<0$, the control constraint $u_i(t)-u_{max}\le 0$ is activated if $t_i^m<\frac{-3v_i^0+\sqrt{9(v_i^0)^2+12u_{max}\cdot(p_i(t_i^m)-p_i(t_i^0))}}{2u_{max}}$.
(ii) For $a_i>0$, control constraint $u_{min}-u_{i}(t)\le 0$ is activated if $t_i^m<\frac{3v_i^0+\sqrt{9(v_i^0)^2+12 \left \| u_{min} \right \|(p_i(t_i^m)-p_i(t_i^0))}}{2 \left \| u_{min} \right \|}$.

\begin{proof}
Based on Lemma 4, the control constraint activation occurs at $t=t_i^0$. For $a_i<0$, we have $ u_i(t_i^0)=a_it_i^0+b_i>u_{max}$.
Without loss of generality, for $t_i^0=0$, we have $b_i=-a_it_i^m>u_{max}$. Substituting $a_i$ from \eqref{eq:a_i}, we obtain $u_{max}\cdot(t_i^m)^2+3v_i^0t_i^m-3(p_i(t_i^m)-p_i(t_i^0))<0$.
Solving the above quadratic term for $t_i^m$, we obtain
\begin{gather}
  (t_i^m+\frac{3v_i^0+\sqrt{9(v_i^0)^2+12u_{max}\cdot(p_i(t_i^m)-p_i(t_i^0))}}{2u_{max}}) \nonumber\\ 
  \,\, \cdot(t_i^m+\frac{3v_i^0-\sqrt{9(v_i^0)^2+12u_{max}\cdot(p_i(t_i^m)-p_i(t_i^0))}}{2u_{max}})<0.\label{eq:48}
\end{gather}
The first term of \eqref{eq:48} is always positive since $\sqrt{9(v_i^0)^2+12u_{max}\cdot(p_i(t_i^m)-p_i(t_i^0))}>0$. Therefore, the second part has to be negative, and the result follows. Following similar steps, the second part of Theorem 3 can be proved for $a_i>0$.
\end{proof}

We have discussed so far the conditions under which the state and control constraints are activated individually. With this information, we can solve the constrained optimization problem and derive corresponding analytical solutions, which we will present in the following section.
However, there is a possibility that the constrained optimal solution may result in additional activation of constraint arcs, i.e., a state constrained optimal solution might cause a control constraint activation. Similarly, a control constrained optimal solution might cause a state constraint activation. In such cases, the optimization problem has to be resolved by piecing all constrained arcs together to derive the corresponding solution. To overcome this recursive procedure, we can identify a priori under which conditions whether other constraint arcs might be activated.

\textbf{Theorem 4.}
(i) The constrained optimal solution corresponding to the state constraint $v_i(t)-v_{max}\le0$ with a junction point at $t=\tau_s$ may result in activating the control constraint $u_i(t)-u_{max}\le0$,  if $  \tau_s < \frac{-3v_i^0+\sqrt{9(v_i^0)^2+12u_{max}\cdot(p_i(\tau_s)-p_i(t_i^0))}}{2u_{max}}$.
(ii) The constrained optimal solution corresponding to state constraint $v_{min}-v_i(t)\le0$ with a junction point at $t=\tau_s$ may result in activating the control constraint $u_{min}-u_i(t)\le0$, if $ \tau_s<\frac{3v_i^0+\sqrt{9(v_i^0)^2+12 \left \| u_{min} \right \|(p_i(\tau_s)-p_i(t_i^0))}}{2 \left \| u_{min} \right \|}$.

\begin{proof}
Suppose that the state constrained ($v_i(t)-v_{max}\le0$) solution has a junction point between the constrained and unconstrained arc at $t=\tau_s$, where $u_i(\tau_s)=0$. Therefore, we focus on the unconstrained arc of the constrained solution with the time horizon $[t_i^0,\tau_s]$ instead of $[t_i^0,t_i^m]$. Following similar procedure as in the proof of Theorem 3, we can construct the condition of the first part of Theorem 4 that indicates whether the control constraint $u_i(t)-u_{max}\le0$ is being activated within the time horizon $t\in [t_i^0,\tau_s]$. The second part of Theorem 4 can be proved following similar arguments.
\end{proof}

\textbf{Theorem 5.}
(i) The constrained optimal solution corresponding to the control constraint $u_i(t)-u_{max}\le0$ with a junction point at $t=\tau_c$ may result in activating the state constraint $v_i(t)-v_{max}\le0$, if $    t_i^m<\tau_c + \frac{3(p_i(t_i^m)-p_i(\tau_c))}{v_i^0+u_{max}\tau_c+2v_{max}}$.
(ii) The constrained optimal solution corresponding to the control constraint $u_{min}-u_i(t)\le0$ with a junction point at $t=\tau_c$ may result in activating the state constraint $v_{min}-v_i(t)\le0$, if $ t_i^m>\tau_c + \frac{3(p_i(t_i^m)-p_i(\tau_c))}{v_i^0+u_{min}\tau_c+2v_{min}}$.

\begin{proof}
Suppose that the control constrained ($u_i(t)-u_{max}\le0$) solution has a junction point between the constrained and unconstrained arc at $t=\tau_c>t_i^0$. Therefore, we focus on the unconstrained arc of the control constrained solution with the time horizon $[\tau_c,t_i^m]$ instead of $[t_i^0,t_i^m]$. Following similar procedure to the proof of Theorem 2, we can construct the condition stated in the first part of Theorem 5 which indicates whether the state constraint $v_i(t)-v_{max}\le0$ is being activated within the time horizon $t\in [\tau_c,t_i^m]$. Similar approach can be employed to prove the second part of Theorem 5.
\end{proof}

\section{Analytical Solution of Constrained Optimization} \label{sec:analytic-solution}
To derive the analytical solution of \eqref{eq:decentral_problem} using Hamiltonian analysis, we adopt the standard procedure used in optimal control problems with control and state constraints \cite{bryson1975applied}, and present a simplified framework. We first start with the result of Theorem 1 to reduce the set of possible constraint activations. Then we evaluate the conditions in Theorem 2 and 3 to identify possible constraint activations. If no activation is identified, we simply derive the unconstrained solution. If any contraint activation is detected by applying the conditions in Theorem 2 and 3, we then use the conditions in Theorem 4 and 5 to identify any additional constraint activations that might rise from the previous case. Once the nature of the current and additional constraint activations is specified, we then piece together the relevant unconstrained and constrained arcs that yield a set of algebraic equations which are solved simultaneously using the boundary conditions of and interior conditions between the arcs.

\textbf{Case 1:}
If only the state constraint $v_i(t)-v_{max}\le 0$ is activated, we have $\mu^{a}%
_{i} = \mu^{b}_{i}=\mu^{d}_{i}=0$. The corresponding necessary condition for optimality, and the Euler-Lagrangian equations \eqref{eq:EL1} for the costates  become $u_i+\lambda_i^v =0$, 
     $\dot{\lambda}_i^p=0$, and
     $\dot{\lambda}_i^v=-\lambda_i^p-\mu_i^c$.
Let us assume that the time $t=\tau_s$ such that $t_i^0<\tau_s<t_i^m$, the unconstrained arc enters the constrained arc. We denote $\tau_s^-$ and $\tau_s^+$ as the immediate left and the right side of $\tau_s$. We have an tangency condition ${N}(t,\textbf{x}_i(t))=v_i(t)-v_{max}=0$, which yields the optimal state $v_i^*(t) = v_{max}$ and control $u_i^*(t) = 0$ on the constraint arc. Considering the jump conditions of the costates and Hamiltonian at the entry point as follows,
\begin{gather}
    \boldsymbol{\lambda}_i(\tau_s^-)-\boldsymbol{\lambda}_i(\tau_s^+)={\eta}_i\frac{\partial {N(t,\textbf{x}_i(t))}}{\partial \textbf{x}_i}\bigg|_{t=\tau_s}, \label{eq:jump-costate}\\
     H_i(\tau_s^+)-H_i(\tau_s^-)={\eta}_i\frac{\partial {N(t,\textbf{x}_i(t))}}{\partial t}\bigg|_{t=\tau_s}, \label{eq:jump-hamiltonian}
\end{gather}
where, ${\eta}_i$ is a constant Lagrange multiplier.
Note that, the state  variables remain continuous, i.e., $\textbf{x}_i(\tau_s^-)=\textbf{x}_i(\tau_s^+)$. From the above equations, both the position costate and the Hamiltonian are continuous at $t=\tau_s$. With $u_i(\tau_s^+)=0$, \eqref{eq:jump-hamiltonian} yields $u_i(\tau_s^-)=u_i(\tau_s^+)$, which implies continuity of control.
Using \eqref{eq:EL1}, \eqref{eq:nece_cond}, state and control continuity, the initial and final conditions \eqref{eq:decentral_problem}, and the terminal condition of co-states at $t=t_i^m$, we can formulate a closed form analytical solution piecing the unconstrained and constrained arcs together at $t=\tau_s$.

\textbf{Case 2:}
If only the control constraint $u_i(t)-u_{max}\le 0$ is activated, we have $\mu^{b}_{i}= \mu^{c}_{i}=\mu^{d}_{i}=0$. The corresponding necessary condition for optimality, and the Euler-Lagrangian equations \eqref{eq:EL1} for the costates  become $u_i+\lambda_i^v+\mu_i^a =0$, 
     $\dot{\lambda}_i^p=0$, and
     $\dot{\lambda}_i^v=-\lambda_i^p.$
Let us assume that at a time $\tau_c>t_i^0$, the unconstrained arc leaves the constrained arc. We denote $\tau_c^-$ and $\tau_c^+$ as the immediate left and the right side of $\tau_c$. For pure control constraint case, we do not have any discontinuity of the costates and Hamiltonian at $t=\tau_c$, i.e., $\boldsymbol{\lambda}_i(\tau_c^-)=\boldsymbol{\lambda}_i(\tau_c^+)$ and $H(\tau_c^+)=H(\tau_c^-)$. Solving these jump conditions, we have $u_i(\tau_c^+)=u_i(\tau_c^-)=u_{max}$ implying the continuity of control input at $t=\tau_c$.
Using \eqref{eq:EL1}, \eqref{eq:nece_cond}, state and control continuity, the initial and final conditions \eqref{eq:decentral_problem}, and the terminal condition of co-states at $t=t_i^m$, we can formulate a similar closed form solution as mentioned in case 1.

\textbf{Case 3:}
If both the control constraint $u_i(t)-u_{max}\le 0$ and state constraint $v_i(t)-v_{max}\le 0$ are activated, we can derive the analytical solution following the procedure described in the previous two cases.

The remaining constraint activation cases for only $v_{min}-v_i(t)\le 0$, only $u_{min}-u_i(t)\le 0$ and their combination can be derived following similar procedure as described above.

\section{Simulation Results}
We validate the analytical solution of the constrained optimization problem through numerical analysis in MATLAB. We select the initial and final position as $p_i(t_i^0)=0\,$ m and $p_i(t_i^m)=200\,$ m, and the initial speed as $v_i(t_i^0)=13.4\,$ m/s. We enforce maximum speed $v_{max}=21\,m/s$ and maximum acceleration $u_{max}=1.4\,$ m/s$^2$. We present only the cases with the state constraint $v_i(t)-v_{max}\le 0$ and control constraint $u_i(t)-u_{max}\le 0$ activation.

\begin{figure}[ht]
\centering
\includegraphics[width=3.5in]{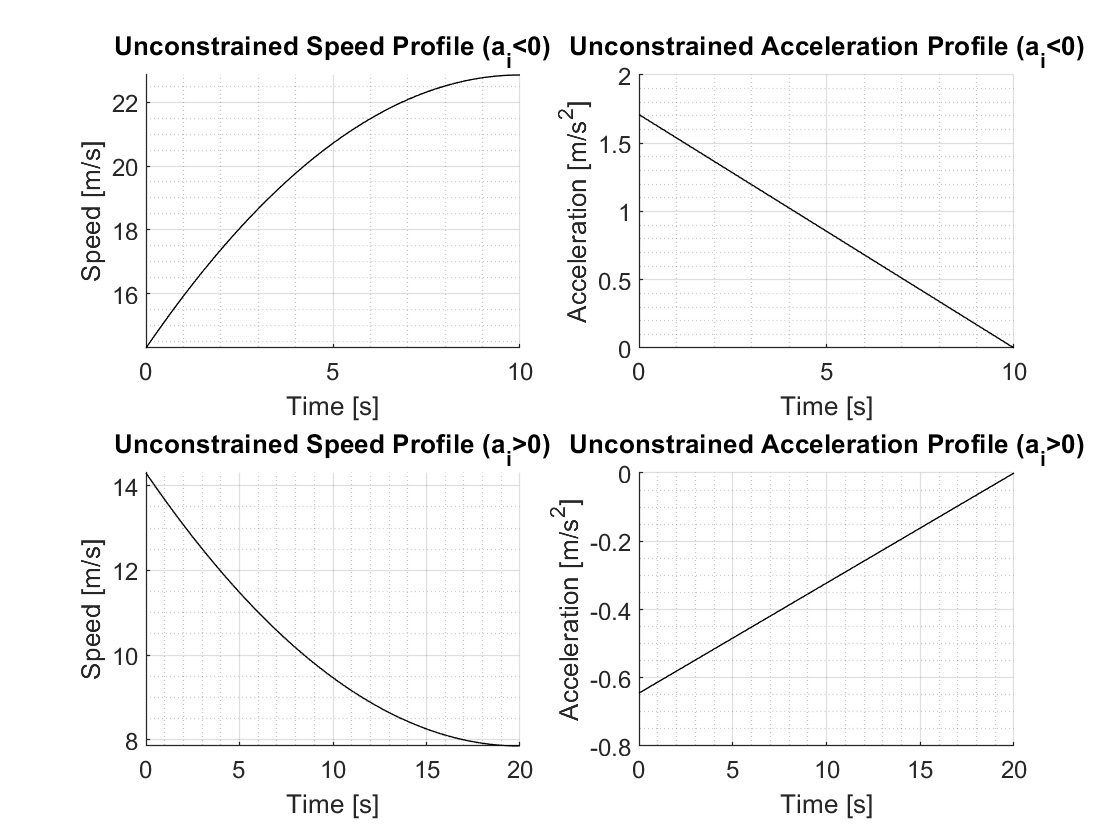} \caption{State and control unconstrained optimal profile for 1) $a_i<0$ (top) and 2) $a_i>0$ (bottom) cases.}%
\label{fig:1}%
\end{figure}

\begin{figure}[ht]
\centering
\includegraphics[width=3.5in]{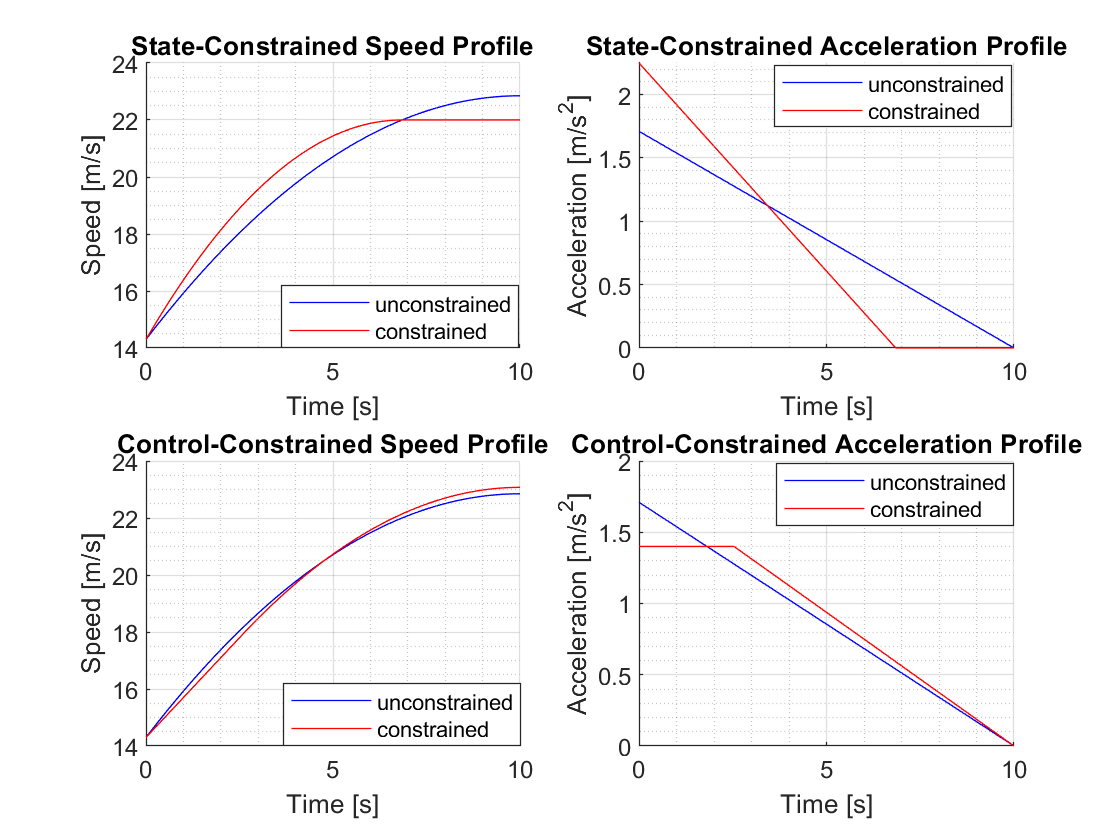} \caption{State constrained optimal speed (top-left) and acceleration (top-right) profile, and control constrained optimal speed (bottom-left) and acceleration (bottom-right) profile.}%
\label{fig:2}%
\end{figure}

\begin{figure}[ht]
\centering
\includegraphics[width=3.5in]{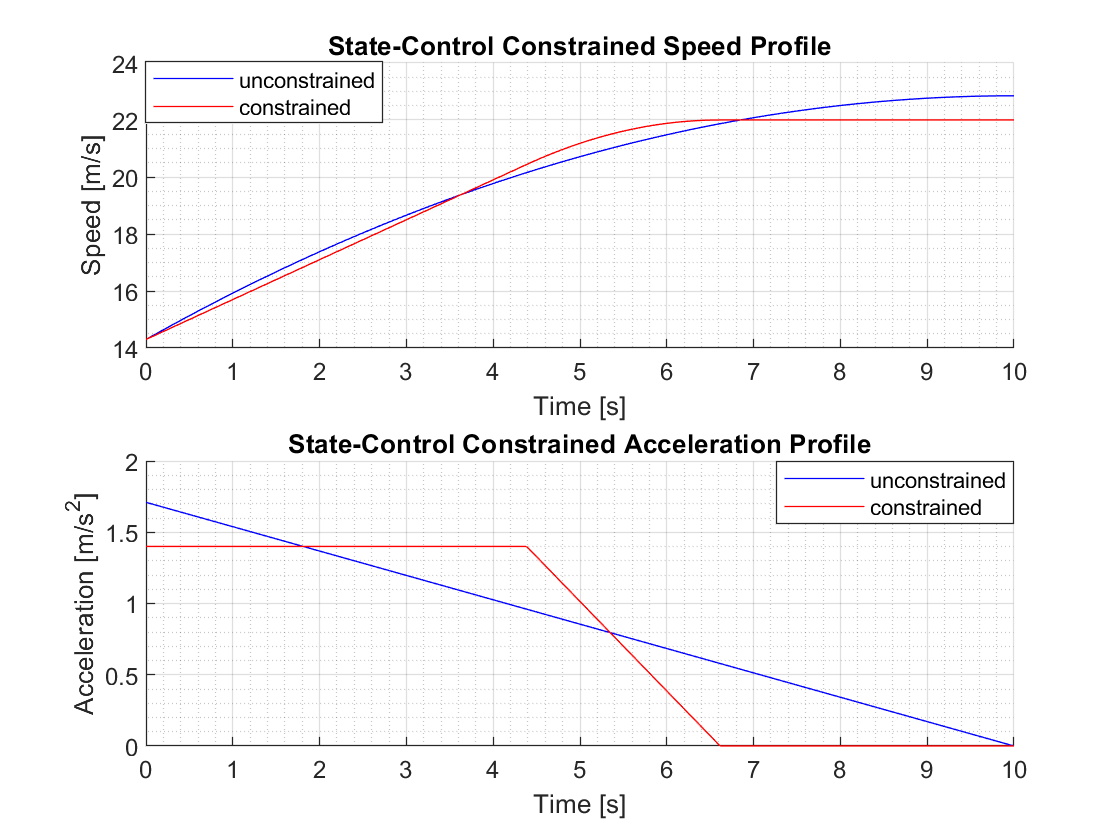} \caption{State and control constrained optimal speed (left) and acceleration (right) profile.}%
\label{fig:3}%
\end{figure}

The unconstrained optimal trajectory for two different merging time $t_i^m$ at $10\,$s and $20\,$s is shown in Fig. \ref{fig:1}. As stated in Lemma 2, we observe two different types of optimal trajectories. This implies that based on the terminal conditions, a CAV may speed up, or slow down optimally to satisfy the boundary conditions. We also observed that both the predefined state and control constraints (for $a_i<0$) are activated in Fig. \ref{fig:1} (top-left and right).
We can save  computational effort for producing the unconstrained results in Fig. \ref{fig:1}, if we know that any of the constraints might be activated. Once we employ Theorem 1 to reduce possible constraint activation set, we use Theorems 2 and 3 to pinpoint the specific constraint activation case. We illustrate the optimal state and control trajectories for only state constrained in Fig. \ref{fig:2} (top) and only control constrained in Fig. \ref{fig:2} (bottom) cases. However, the results in Fig. \ref{fig:2} have to be recalculated if additional constraint activation is encountered. We can avoid the unnecessary computational effort to produce the intermediate result in Fig. \ref{fig:2} based on the following observation. 
Note that, in the state constrained solution, the control constraint is activated (Fig. \ref{fig:2}, top-right) which was non-existent before, as discussed in theorem 4. Similarly, we observe in Fig. \ref{fig:2} (bottom-right) that the control constrained optimal trajectory creates a possibility of additional state constraint activation (Fig. \ref{fig:2}, bottom-left) due to the increased speed, as discussed in Theorem 5. Hence, we need to take the additional constraint activation into account from Theorems 4 and 5, and enforce both state and control constraints if needed. The unconstrained and fully constrained state and control constrained trajectories are shown in Fig. \ref{fig:3}. The four-point boundary value problem is solved in this case, and two constrained and one unconstrained arcs are pieced together to provide the optimal solution.

\section{Concluding Remarks}
In this paper, we presented a theoretical analysis of the decentralized optimal control framework for coordinating CAVs that was presented earlier \cite{Malikopoulos2017}. We
provided the conditions under which the state and control constraints are activated along with the corresponding solutions. These conditions can provide insight about when any of the state and control constrained arc is activated, and thus can be useful for real-time implementation.
Ongoing work includes the analysis of other cases, including safety constraints, along with the conditions that a solution always exists. A potential direction of future research includes the enhancement of this analysis by considering different penetration rates of CAVs and the implications of error or delays in the communication.

\bibliographystyle{IEEEtran}
\bibliography{acc_pt_vd_ref}

%

%

\end{document}